\documentclass[12pt,notitlepage,twoside,a4paper]{amsart}
 \usepackage{amsfonts} 
\usepackage[latin1]{inputenc}
\usepackage[cyr]{aeguill}
\usepackage{xspace}
\usepackage[polutonikogreek,english]{babel}

\usepackage{amsmath,amssymb,enumerate}

\usepackage{epsfig,fancyhdr,color}%,showkeys,amsmidx

\usepackage{amssymb}
\usepackage{amsmath,amsthm}
\usepackage{latexsym}
\usepackage{amscd}
\usepackage{psfrag}
\usepackage{graphicx}
\usepackage[latin1]{inputenc}
\usepackage[all]{xy}
\usepackage[mathcal]{eucal}

\definecolor{NoteColor}{rgb}{1,0,0}

\renewcommand{\textsc}{\textcolor{red}}

\newtheorem*{theorem 1}{\rm\bf Proposition 1}
\newtheorem*{theorem 2}{\rm\bf Proposition 2}

\theoremstyle{definition}

\theoremstyle{remark}

\def\interieur#1{\mathord{\mathop{\kern 0pt #1}\limits^\circ}}

\title[Clairaut, Euler and the figure of the Earth]{Clairaut, Euler and the figure of the Earth}

\author{Athanase Papadopoulos}
\address{Athanase Papadopoulos,  Universit{\'e} de Strasbourg and CNRS,
7 rue Ren\'e Descartes,
 67084 Strasbourg Cedex, France}
\makeindex
 
\date{\today}

\begin{document}

  \begin{abstract}

      The sphericity of the form of the Earth was questioned around the year 1687, primarily, by Isaac Newton who deduced from his theory of universal gravitation that the Earth has the form of a spheroid flattened at the poles and elongated at the equator. In France, somepreeminent geographers were not convinced by Newton's arguments, and about the same period, based on empirical measurements, they emitted another theory, claiming that on the contrary, the Earth has the form of a spheroid flattened at the equator and elongated at the poles.
      To find the real figure of the Earth became  one of the major questions that were investigated by geographers, astronomers, mathematicians and other scientists in the eighteenth century, and the work done around this question had an impact on the development of all these fields.
            
  In this paper, we review the work of the eighteenth-century French mathematician, astronomer and geographer Alexis-Claude Clairaut related to the question of the figure of the Earth. We report on the relation between this work and that of Leonhard Euler.  At the same time, we comment on the impact of the question of the figure of the Earth on mathematics, astronomy and hydrostatics. Finally, we review some later mathematical developments that are due to various authors that were motivated by this question. It is interesting to see how a question on geography had such an impact on the theoretical sciences.
  
  The final version of this paper will appear in  \emph{Ga\d{n}ita Bh\=ar\=at\=\i }, (Indian Mathematics) the Bulletin of the Indian Society for History of Mathematics.

     

     \end{abstract}
          
  \maketitle
  
  \tableofcontents

  \section{Introduction}

  The question of the figure of the Earth has surely haunted man since early times.  Written documents attest that since the 6th century BCE, astronomers, philosophers and mathematicians of ancient Greece believed that the Earth is spherical. There are probably  traces of similar theories in other cultures. One can imagine the contentment of the astronomers who realized for the first time the non-flatness of the Earth, probably deducing this fact from their observations of the variations in the positions of the stars in the heavenly sky, when they traveled north or south.\footnote{We are not talking here of the the local irregularities caused by the existence of mountains, valleys, seas, and other differences in altitude, but about the shape of the Earth when these differences are considered as negligible compared to its size.} My aim in this paper is to review the question of the figure of the Earth as it was discussed  in the eighteenth century, and to talk about its impact on mathematics and  other science. I will start by recalling why this period was a turning point for this subject.
  
  Around the year 1670, the sphericity of the Earth was questioned. Geographers and physicists, especially in England and France, started to believe that the Earth is not spherical but spheroidal. A spheroid is a figure obtained by the rotation of an ellipse around one of its axes.  Two conflicting theories arose, the one led by Newton, who deduced from his theory of universal gravitation that the Earth is a spheroid flattened at the poles and elongated at the equator, while the other one, supported by the famous astronomers of the Cassini family  and their collaborators in France, concluded, based essentially on experimental data, that the Earth, on the contrary, is an ellipsoid elongated at  the poles and flattened at the equator.  Newton's conclusion was based on his conviction that the Earth was originally constituted of a fluid having an overall spherical shape, and that it reached its present spheroidal form as a consequence of its rotation around its axis and under the effect of the difference in the gravitational field at the poles and the equator. He even predicted in his \emph{Philosophiae Naturalis Principia Mathematica} (1687) a precise value for the flattening, namely, 1/230, meaning that if $a$ is the major axis and $b$ the minor axis of the elliptical section of the Earth, then,  $\frac{a-b}{a}=\frac{1}{230}$ (\emph{Principia}, Book III, Propositions XVIII-XX). 
  
   These two conflicting theories led to a debate that involved some of the most preeminent scientists of Europe and that ended only in 1740 where, after a series of measurements of lengths of degrees of meridians that were conducted near the poles and near the equator, it was agreed that the Earth has the form of a spheroid which is flattened at the poles and elongated at the equator. 
   
   The arguments that were brought to resolve this  question had enormous repercussions on the development of science. Leonhard Euler took part in this  controversy, and this had an impact on his work in mathematics, physics, astronomy and geography. We recall incidentally that Euler, besides being a mathematician (admittedly,  the most preeminent mathematician of his times, and maybe of all times), was also a geographer. He was the main editor of two of the most  important atlases of the eighteenth century, published by the Russian and the Prussian Academies of Sciences, and he wrote several memoirs on the problem of drawing geographical maps and other geographical questions. He was also thoroughly involved in astronomy, and the discovery that the figure of the Earth was not spherical had consequences on his work in this field. He also worked on fluid mechanics, and the various theories emitted on the way the Earth attained its present shape led him to study the mathematical evolution theory of the fluid that constitutes the internal matter of this planet and other questions related to the actual figure of the Earth.
   
    From the purely mathematical point of view, it suffices to note that Euler, motivated by the spheroidal shape of the Earth, and  after having written several memoirs on spherical trigonometry (see the survey in \cite{Euler-GB}),  published a memoir on the trigonometry of the spheroid, which is the earliest known work dedicated entirely to this subject. We shall review some results of this paper in \S \ref{s:Euler} below.

Euler corresponded with several mathematicians and  scientists on the question of the figure of the Earth. His most important correspondent regarding this subject was probably the famous French mathematician, geographer and physicist Pierre-Louis de Maupertuis who, contrary to the opinion of his colleagues at the French Royal Academy of Sciences, was a supporter of the theory claiming that the Earth is fattened at the poles. Maupertuis was the director of the Prussian Academy of Sciences during the period where Euler worked there. We comment at length on Maupertuis' work on the question of the figure of the Earth in relation with that of Euler in the book \cite{CP}. In the present paper, I would like to survey some of Euler's work on this subject in relation with that of Alexis-Claude Clairaut, another major mathematician of the eighteenth century who had a sustained correspondence with Euler. Most of all, I would like to emphasize the influence on the research around the figure of the Earth on the mathematical works of both men.

The plan of the rest of this paper is the following.
I will start with a short review of the question of the Figure of the Earth in Greek antiquity (\S \ref{s:Question}). Then I will pass to the same question in the seventeenth and eighteenth centuries (\S \ref{s:seven}).
 I will then give a glimpse of the life of Clairaut (\S \ref{s:vita}), followed by a quotation from the eighteenth century encyclopedia of Diderot and d'Alembert on Clairaut's work on this topic (\S \ref{s:enc}). I will then report  on the work of Euler on the same subject (\S \ref{s:Euler}), and then review part of the correspondence between the two men on this question (\S \ref{s:correspondence-C}). The fact that the Earth was thought to be spheroidal acted as a motivation for the study of  the  geometry of the spheroid and other surfaces of revolution by later mathematicians,  and I will mention some of the work done on this subject (\S \ref{s:works}).

 \section{On the question of the figure of the Earth in Greek Antiquity}\label{s:Question}
  
  Let me start with some quotes from Greek Antiquity.  Aristotle (384-322 BCE), in \emph{On the heavens} \cite{Aristotle-Heavens} Book II,  Chapter 13, mentions several philosophers among his predecessors who  considered that the Earth is spherical. The names he quotes include Anaximenes (c. 585-525 BCE), known for his reflections on the concept of infinite, Anaxagoras (c. 500-428 BCE), who presumably was the founder of the first philosophical school in Athens, and Democritus (460-370 BCE), who formulated an atomic theory of the universe. 
  Aristotle writes in the same passage that there was an emulation among Greek philosophers to establish who was the first to declare the sphericity of the Earth.
    A similar dispute is reported on by the Greek biographer Diogenes Laërtius  in his book \emph{The lives of the philosophers}  (3d. c. CE) \cite[Book VIII, 1.26]{Diogenes}:  ``We are told that [Pythagoras] was the first to call the heaven the universe and the Earth spherical though Theophrastus says it was Parmenides, and Zeno that it was Hesiod." In the same chapter, Diogenes  quotes a passage from the Roman historian  Alexander Polyhistor (1st c. BCE) who writes in the \emph{ Successions of Philosophers} that the Pythagoreans knew that the Earth is spherical and inhabited round about, explaining in this way the succession of light and darkness, of hot and cold, and of dry and moist. We do not dwell too much on this dispute; our aim here is just to recall that the issue was debated in early times.

   Aristotle, in Chapter 14 of Book II of his treatise \emph{On the heavens}  \cite{Aristotle-Heavens} is concerned with the Earth: its position in the universe, its motion, its form, and its magnitude.  It is there that we find the evidence that the Stagirite gave for the sphericity of the Earth.  
    
    As always, Aristotle first quotes the opinions of his predecessors. Regarding position and form, some of them, he says, considered the Earth as a star among others, while others thought that is at the center of the universe, adding that it has a rotational motion around a central axis.  It is interesting that Aristotle finds that the sphericity is incompatible with a permanent uniform rotary motion. He then gives the reasons for which he considers that the Earth is motionless and at the centre of the cosmos. It would be interesting to enter into the details of his (mathematical and astronomical) arguments, but this is not subject of our paper.  We will neither talk about Aristotle's opinion on the Earth's magnitude.
         We come to the sphericity. Let us first note that
Aristotle's arguments for this sphericity were still not outdated in the eighteenth century. D'Alembert wrote an article titled 
\emph{Figure de la terre} (Figure of the Earth) in the  sixth volume of the   \emph{Encyclopédie}, published in 1756, in which he writes: 
``[Aristotle]  establishes and proves the roundness of the Earth in his second book of \emph{De Caelo}, chap. xiv., by very solid arguments, which are similar to those we are about to give." 
In \S \ref{s:enc}, I will quote more from d'Alembert's \emph{Encyclopédie} article.

 Aristotle's argument goes as follows : Every portion of the Earth has a weight, which attracts it to the center.   This is the famous ``natural downward movement" of objects, of which he often speaks in his writings. He then says that since the weight of a certain mass is the same everywhere on the surface of the Earth, there is a unique center towards which heavy objects are attracted.  He elaborates on this, adding that during its generation  the Earth was already spherical. He writes that this figure is compatible with the observations of the eclipses of the moon, since during such an eclipse, we see the figure of the Earth on the moon.

   A later major scientist from Greek Antiquity, the
famous geographer Strabo (c. 60 BCE - 20 CE), also gave some physical support for the roundness of the Earth. It is interesting to read passages from his work in which first principles ruling the physical world   combined with practical observations imply this roundness. We read in his \emph{Geography} 1.1.20 \cite[vol. 1, p. 41]{Strabo}:  
\begin{quote}\smaller

[\ldots] I must take for granted that the
universe is sphere-shaped, and also that the Earth's
surface is sphere-shaped, and, what is more, I must
take for granted the law that is prior to these two
principles, namely that the bodies tend toward the
centre;  and I need only indicate, in a brief and
summary way, whether a proposition comes---if it
really does---within the range of sense-perception or
of intuitive knowledge.
 Take, for example, the proposition that 
the Earth is sphere-shaped: whereas the suggestion of this 
proposition comes to us mediately from the law that 
bodies tend toward the centre and that each body inclines toward its own centre of gravity, the suggestion 
comes immediately from the phenomena observed at 
sea and in the heavens [\ldots]  our sense-perception and 
also our intuition can bear testimony in the latter 
case. For instance, it is obviously the curvature of 
the sea that prevents sailors from seeing distant lights 
that are placed on a level with their eyes. At any 
rate, if the lights are elevated above the level of the 
eyes, they become visible, even though they be at a greater distance from the eyes; and similarly if the 
eyes themselves are elevated, they see what was 
before invisible. 
\end{quote}

Ptolemy (c. 100-168 CE), in the \emph{Almagest} (Book I, 3-8), \cite{Halma-Ptolemy} follows Aristotle in considering that the Earth is spherical. At  the same time, he mentions other possibilities (plane, polygonal, convex, concave, cylindrical), shwing how astronomical evidence rules them out.  

Finally, let us mention the 4th-century mathematician Pappus of Alexandria, who in his \emph{Mathematical collection} takes advantage of the roundness of the Earth in order to introduce a series of mathematical propositions. He declares in \cite[t. 1, p. 272]{Pappus} that the philosophers (by which he means the mathematicians) rightly said that the first among the gods gave to the Earth a spherical shape, the most beautiful shape, adding that the sphere, among the  surfaces having the same area, encloses the largest volume. He says that the philosophers made this statement without proving it, and he gives a series of propositions in which he proves this statement for the restricted class of polyhedral surfaces. The corresponding general result for surfaces, known as one form of the so-called isoperimetric inequality, was proved only  in the nineteenth century.

 \section{The question of the figure of the Earth in the seventeenth and eighteenth centuries}\label{s:seven}
The sphericity of the Earth started to be challenged during the last three decades of the seventeenth century. In 1669, 
Louis XIV (the ``Sun King") asked the French geographers and astronomers,  led by Jean Picard who was joined later by Jean-Dominique Cassini
  and then by his son Jacques,\footnote{The Cassini family, of Italian origin, included an important number mathematicians, astronomers  and geographers who lived in France  in the seventeenth and eighteenth centuries.  The first one among them is Jean-Dominique Cassini (1625-1712), who was initially a well-known astronomer, mathematician and hydraulic engineer at the University of Bologna. He was called to Paris in 1669 by Jean-Baptiste Colbert, who was the most important secretary of state of Louis XIV, in charge of marine and of finance. Colbert proposed Cassini to join the newly established Royal Academy of Sciences and he offered him the post of director of the Paris Observatory. This was at a time where France was in need of precise measurements of longitudes, latitudes and distances and of scientists capable of drawing accurate maps of all the parts of the kingdom. Cassini arrived indeed with a new method for determining  longitudes. 
He became later known as Cassini I. His son Jacques (1677-1756) became also a respected astronomer and geodesist, and he was known under the name Cassini II. He was followed by César-François Cassini (1714-1784), son of Jacques, also known as Cassini III, and Jean-Dominique Cassini de Thury (1748-1845), son of César-François, also known as Cassini IV. There were other scientists in this family. The so-called Cassini map of France, also known as the Academy map, is considered as the first precise map of the Kingdom of France as a whole and it was drawn using a geodesic triangulation whose constitution took more than sixty years. The drawing of the final map is mainly due to Cassini III and Cassini IV, but the work behind it is due to the four generations of Cassinis.}  to conduct a project of land surveying in France based on the method of triangulations and using astronomical observations. The goal was  to measure the length of the meridian arc between the two cities of Amiens and Paris, with the purpose of drawing precise maps. The scientists concluded from the measurements they conducted that the Earth is not spherical, but spheroidal, that is,  it is obtained by the rotation of an ellipse (and not a circle) around an axis.  In 1718,  Jacques Cassini made an address at the Royal Academy of Sciences in which he declared that the measurements that were led by his father to measure the length of a degree of meridian, conducted at different altitudes in the Northern hemisphere, showed that this length increases as  one moves North, which means that the Earth is elongated at the poles. Thus, he concluded that the Earth is a prolate spheroid, that is, a spheroid flattened at the equator and elongated at the poles. 
  
A theory stating the contrary, namely, that the Earth is an oblate spheroid, that is, a spheroid
elongated at the equator and flattened at the poles,
 was emitted by Isaac Newton in the 1687 edition of his \emph{Principia}  (Book III, Propositions XVIII-XX). In other words, Newton considered that the rotation axis of the Earth, is smaller than the other axis (usually called the \emph{diameter} of the Earth). Furthermore, he gave a precise estimate of the flattening of the Earth, namely,  he declared that if $a$ denotes the major axis and $b$ the minor axis of the elliptical section of the Earth, then $\frac{a-b}{a}=\frac{1}{230}$. His theory was based on the assumption that the Earth was originally a fluid having a spherical shape, and that it acquired its spheroidal shape gradually, under the effect of the mutual attraction force exerted between the various parts of the Earth and under the effect of  the latter's rotation around its axis.

  Back in 1672, that is, before Newton emitted his theory, the French astronomer Jean Richer observed that the length of a  pendulum performing one beat per second is longer in Paris than in Cayenne (a place in French Guiana, near the equator).  The only explanation available  for this difference in length was that the value of gravitational force is smaller in Cayenne than in Paris (after the variation of the length of the pendulum due to temperature has been taken into account).   This would mean that the factor representing the acceleration of gravity  which appears in the formula for the period of the pendulum was not  the same at each place: it depends on the latitude of the location.    As a matter of fact, Richer was sent to Cayenne by Jean-Dominique Cassini who stayed in Paris, so that they could both and simultaneously observe the planet Mars and, by computing its parallax, they could give an estimate of the distance from Mars to the Earth. Newton considered that Richer's results support his theory. 
    
  Christiaan Huygens, who was certainly the most preeminent seventeenth-century mathematician in Paris, was of the opinion of Newton, concerning the figure of the Earth, although he gave a different explanation for this phenomenon.  In fact, Huygens considered that the particles constituting the Earth were submitted at the same time to the force of attraction and to the centrifugal force due to the rotation of the Earth around its axis. He concluded that, in order to reach in equilibrium under these two forces, a flattening  of the Earth towards the poles was necessary. He also gave an estimate of this flattening, which was different from the estimate that Newton gave. The latter, although he made a different reasoning, considered that Huygens' theory on the shape of the Earth was a confirmation of his own theory.

In 1736-1737, the French Academy of Sciences, together with the ministries of Marine and of Finances, organized an expedition to Swedish Lapland whose objective was to measure the length of a degree of meridian in these regions which are close to the North Pole. Comparing the result with the length of a degree of a meridian at some known place (like Paris) would give information on whether the Earth is flattened or not at the pole. The leader of the expedition was Pierre-Louis de Maupertuis, who was a member of the Royal Academy of Sciences and who was at the origin of the idea of this expedition. Maupertuis was one of the main supporters in France of the theory saying that the Earth is oblate. More generally, he was a supporter of Newton's ideas, against the majority of the other French academicians.  We mention incidentally that in those times, there was a competition in France between Newton's and Descartes' ideas on  subjects such as matter and attraction, and the dispute regarding the figure of the Earth was part of the debate between the supporters of the two theories, a debate that sometimes took the form of a conflict. Voltaire, at several places of his pamphlets and his literary essays, mentions this dispute. He writes in a letter to M. de Formont, dated December 23, 1737: ``The spirits are in Paris in a small civil war; the Jansenists attack the Jesuits, the Cassinists rise against Maupertuis, and do not want the Earth to be flat at the poles." \cite[p. 184]{Voltaire}

Another expedition headed to Peru, the year before the one to Lapland, to make similar measures near the equator. The Peru expedition lasted 9 years (1735 to 1744).

Maupertuis became a few years later the director of the Prussian Academy of Sciences where Euler was working, and the two men had very close relations. Clairaut was one of the youngest scientists that took part in the Lapland expedition. We shall talk about Clairaut and Euler in the next sections.

  \section{{\it Vita} of Clairaut}\label{s:vita}
  
  The life of Clairaut is recounted by Jean-Paul Grandjean de Fouchy, an astronomer who was, during 32 years, the perpetual secretary of the Royal Academy of Sciences of Paris, and who, as such, was in charge of writing the obituaries of the deceased members of that Academy. The information on Clairaut's life contained in the present article is collected from Grandjean de Fouchy's obituary \cite{Grandjean-Clairaut} and from a few other sources including the  biographical essay by Pierre Brunet \cite{Brunet}. We have also extracted some information from Clairaut's correspondence with Euler.

Alexis-Claude   Clairaut  (1713-1765) was the second of twenty-one siblings. His father, Jean-Baptiste Clairaut, was a mathematician who taught  in Paris and who was a member of the Royal Academy of Berlin and of famous other learned societies. Alexis-Claude was a precocious and exceptionally talented child.\footnote{On the contrary, Leonhard Euler, who is the other main figure in this paper, is considered to have had a normal child evolution.} He received all his education at home. From his biographers, we learn that Alexis-Claude's father taught him the letters of the alphabet using, as support, the  figures of an edition of Euclid's \emph{Elements}, together with their captions.  The child soon wanted to understand the meaning of the words he was spelling, which was obviously his father's wish. In this way, he became familiar with the geometry of triangles very early in his life. Jean-Baptiste Clairaut used another ingenuity for the young Alexis-Claude to learn Latin, putting at his disposal a collection of books on war machines written in that language, containing many illustrations. This instilled in the child, who was attracted by the drawings, the desire of reading and understanding the text. At age ten, Clairaut was able to read the \emph{Conic sections}  of the Marquis de l'Hôpital, and soon later, books on infinitesimal methods, differential calculus and integration by the same author.

  Jean-Baptiste Clairaut soon became concerned about his child's health as the latter spent his whole days reading and studying mathematics, and he tried to reduce his enthusiasm for this subject. The young Alexis continued to work secretly, together with his younger brother,\footnote{Clairaut's younger brother who, supposedly, was talented like him, died before age seventeen.}   while the rest of the family was asleep, despite the fact that his father, as Grandjean de Fouchy says, ``had seriously prescribed this studious debauchery."
 At the same time, Jean-Baptiste Clairaut did not wish his son to lose the enthusiasm for learning that he had acquired, and he wanted the Academy  of Sciences to be aware of the work of this exceptionally gifted child. 
 In 1726, he presented to the Academy an original work on plane curves done by his 
a thirteen years old son.  The academicians confirmed that the results were interesting but they first received them with great suspicion; indeed, no one at the Academy could believe that these results were due to a 
a thirteen years old child. It was only after a session of questions and answers   that the boy received the praise he deserved.
Clairaut's memoir was published several years later. An English translation of this memoir appeared recently, see  \cite{MAA}.

In 1729,  Clairaut, who was sixteen,  wrote a memoir titled \emph{Recherches sur les courbes à double courbure} (Researches on curves with double curvature) \cite{C0} in which he studies space curves. Space curves at that time, were generally considered as curves drawn on the surface of a solid. An example of a space curve is a curve obtained by making a compass turn on a cylinder or on another surface. The subject was relatively new.\footnote{The name ``courbes à double courbure", used by Clairaut, had already been used a few years before by Henri Pitot (1695-1771) in a memoir
titled \emph{Sur la quadrature de la moitié d'une courbe qui est la compagne des arcs, appelée la
compagne de la cycloïde} (On the quadrature of half of a curve which is the companion of arcs,
called the companion of the cycloid), presented to the Académie Royale in 1724 and published two
years later \cite{Pitot}.  Pitot considered a spiral on a cylinder. He writes about it:  ``The Ancients called this curve
spiral or helix, because its construction on the cylinder follows the same analogy as the construction
of the ordinary spiral on a plane, but it is very different from the ordinary spiral, being one of these
curves with double curvature or a line which one can conceive as traced out
on the curved surface of a solid."}
 Clairaut  considered that he was the first to investigate such curves. He introduced a technique for studying them  by considering the projection of such a  curve on two planes making a right solid angle. He declared that 
Descartes said that to investigate such a curve, one has to project it onto two
perpendicular planes, and transform it into curves contained in these planes. He considered that Descartes, who planned such a study, did not carry it out.

 Clairaut worked so intensively on that memoir that he became sick, and it took him two years to recover. Grandjean de Fouchy writes that in that paper, Clairaut developed a principle  which ``opened to the geometers a new stone-pit in which no one had been able or wished to engage until then"  \cite[p. 149]{Grandjean-Clairaut}. Clairaut's memoir was published in 1731 by the Royal Academy, with a certificate, signed by Fontenelle, perpetual secretary of the Academy, testifying that this institution has taken every precaution to ensure that the author was barely 16 years old when he submitted this work. The same year (1731),  Clairaut was elected at the Academy as an adjunct in mechanics. A permission from the King, Louis XV, was needed, since the lowest age permitted was 20. Clairaut was only 18.

Soon after Clairaut  was elected at the Academy, he presented  two geometrical memoirs, the one titled \emph{Nouvelle manière de trouver les formules des centres de gravité} (A new manner to find the formulae for centers of gravity) \cite{C011} and the other one, \emph{Sur les courbes que l'on forme en coupant une surface courbe quelconque par un plan donné de position} (On the curves formed by cutting an arbitrary surface by a given plane of position). The latter may be considered as a sequel to his 1729 memoir on space curves. In this new memoir, Clairaut studies surfaces by intersecting them with planes. The two memoirs were published two years later in the Academy's Memoirs series, \cite{C01, Sur-1733}.

Between September and November 1734, Clairaut stayed in Basel where he followed the teaching on differential calculus given  by  Johann I Bernoulli, who was one of the most renowned mathematicians of his times. Bernoulli had been the teacher of Leonhard Euler, and was certainly  the best specialist in infinitesimal calculus of his times.   Clairaut went to Basel with Maupertuis, who had already  strong relations with the Bernoulli family. He had made a stay there  in 1729-1730, where he also followed lessons by Johann I Bernoulli, and after that he kept a regular correspondence with Johann II Bernoulli, one the the children of Johann I, who was staying with him in Basel. (Nikolaus  and Daniel Bernoulli,  the other two children of Johann I, had moved to Saint Petersburg in 1726.) Maupertuis, when he decided to return to Basel in 1734, proposed to Clairaut  to accompany him, and the latter did not hesitate to do so.  Maupertuis was 15 years older than Clairaut, and he was his friend and colleague at the Royal Academy.

 The contact between Clairaut and Maupertuis is  important for the subject of our paper, because 
 Maupertuis was already thoroughly involved in the question of the figure of the Earth.  Grandjean de Fouchy writes that at the time Clairaut returned from Basel, the Academy was so busy with this question  that  it was natural that Clairaut got also involved in it. 
 
 Soon after they returned from Basel, Clairaut and  Maupertuis went for a retreat in  Mount Valérien, a hill near Paris, known for its calmness, with the aim of discussing the issue of the figure of the Earth.\footnote{Mount Valérien is a hill situated West of Paris which has been  a hermitage until the French revolution. Later, it remained a place of devotion. Jean-Jacques Rousseau used to go there for his meditations.} The place was also propice for mathematical work.  It seems that it is during this retreat that Maupertuis formed the project of an  expedition to Lapland, on the Arctic circle, whose aim was to make precise computations of the length of a degree of a meridian near the North Pole, in order to confirm the fact that the Earth has the form of a spheroid which is flattened at the poles.

In 1735, upon the recommendation of their common friend Daniel Bernoulli, Euler tried to hire Clairaut at the Saint Petersburg Academy of Sciences, but his wish was not fulfilled,  see \cite[p. 2]{Opera-4-A-7}.

Clairaut  was one of the important members of the 1736--1737 Lapland expedition which was directed by Maupertuis.   In fact, he was one of the four members of the Royal Academy of Sciences that were part of this expedition, the others being Maupertuis himself, the mathematician and astronomer Louis Camus, and the astronomer and geodesist Pierre Charles Le Monnier, who was a young adjunct at the Academy. The expedition lasted 16 months. 
The Swedish astronomer Anders
Celsius, who was a professor at the University of Uppsala, also joined the expedition.

 In 1740, Clairaut wrote his first letter to Euler.    This was the beginning of an extensive correspondence which lasted  until 1764, in which the two men discussed a variety of topics including the integration of special functions, the differential geometry of curves, isoperimetry problems, differential equations,  optical instruments, hydrodynamics  and celestial mechanics, in particular the motion of the moon.

In his first letter to Euler, dated September 17, 1740,  \cite[p. 68]{Opera-4-A-7}, Clairaut writes:\footnote{In this paper, unless otherwise indicated, the translations from the French are mine.}
\begin{quote}\small
I have wished since an infinite time to be in correspondence with you, but the fear of appearing too daring and that of diverting you from your occupations have prevented me from taking the liberty of writing to you until now. Although these reasons still remain, I envy so much the pleasure that several learned friends of mine\footnote{Daniel Bernoulli and Maupertuis were surely among these friends.} that I can no longer prevent myself from granting the same grace to myself as to them by honoring myself from time to time with your letters. If the little geometry that has appeared from me has reached you and has merited enough of your attention so as to obtain for me the grace I ask for, I will be comforted to have produced so little so far.
 \end{quote}

In his response to Clairaut, dated October 19, 1740  \cite[p. 71]{Opera-4-A-7}, Euler writes:  
\begin{quote}\small
Admiring for a long time already your penetrating intelligence, I burned with the desire to maintain a correspondence with you, being certain that I would make great progress with the help of your profound reflections.  It is with great pleasure, Sir, that I have read your remarkable works: not only those which have been published separately or inserted in the \emph{Mémoires} of your Academy, but also, very recently, your articles of a rare quality on the figure of the Earth with which you have enriched the English.\footnote{Euler refers to articles in English that Clairaut published in the Transactions of the Royal Society.}
\end{quote}

The question of the figure of the Earth is naturally one of 
 the subjects the  two scientists addressed in their correspondence, and we shall review part of this correspondence on this subject in \S \ref{s:correspondence-C}. Of course,
geometry was also one of their major common interests. Euler mentions Clairaut at the beginning of  the first chapter of the second volume of his important treatise \emph{Introductio in analysin infinitorum} (1748), in which he studies curvature of space curves: ``Curves of this kind have two kinds of curvature  which has been beautifully discussed by the brilliant geometer Clairaut".

In 1741, Clairaut published a geometrical treatise, \emph{\'Eléments de géométrie} (Elements of geometry) \cite{Elements-1741}. He writes in the Preface (p. vii): ``Land measurement seemed to me to be the most appropriate thing to give rise to the first propositions of Geometry; and it is indeed the origin of this science, since \emph{Geometry} means \emph{measurement of the Earth}." The first \emph{Preliminary notion} of his treatise reads: \emph{The way we measure an arbitrary length.} In the same preface, Clairaut gives interesting thoughts on the difference between his \emph{Elements} and those of Euclid. We quote a passage which will give the reader an idea about his style and his way of thinking. He writes (p. ix-xx):
\begin{quote}\small
I may be accused, at some places in these Elements, of relying too much on the testimony of the eyes, and of not attaching enough importance to the rigorous exactitude of the demonstrations. I beg those who might make such a reproach to observe that I only pass briefly over propositions whose truth can be discovered as long as one pays attention to them. I do this, especially in the beginning, when such propositions are more often encountered, because I have noticed that those who were inclined to Geometry enjoyed exercising their minds a little, and that on the contrary they were repulsed when they were burdened with demonstrations that were, so to speak, useless.
 
When Euclid takes the trouble of demonstrating that two intersecting circles do not have the same center, that a triangle enclosed in a second one has the sum of its sides smaller than that of the sides of the triangle in which it is enclosed, no one will be surprised.  This geometer had to convince obstinate sophists, who prided themselves on refusing to accept the most obvious truths: it was therefore necessary that Geometry should have, like Logic, the help of formal reasoning to close the mouth of the dispute. But things have changed. Any reasoning that falls on what common sense alone decides in advance is nowadays pure waste, and is only suitable to obfuscate the truth and to disgust the readers.

\end{quote}

In 1743,  Clairaut published his major work  \emph{Th\'eorie de la figure de la terre, tir\'ee des principes de l'hydrostatique} (Theory of the figure of the Earth, extracted from the principles of hydrostatics) in which, following Newton, he assumes that the Earth was initially a fluid, studying the equilibrium state of that fluid submitted to the law of universal attraction and to the centrifugal force due to the Earth rotation around an axis. He treated first the case of a fluid of fixed density, and then that of variable density, constant on strata. 
In this book, Clairaut formulated a mathematical theory of hydrodynamics that had an enormous impact on science; we shall say more about this in \S \ref{s:works} below.
The work was presented to the Royal Academy at several sessions in 1741 and 1742.

In the same year (1743), Clairaut  started a major work on the motion of the moon, again motivated by Newton's  theories. This motion involves the simultaneous mutual attraction of the sun,  Earth and moon, following Newton's law of universal gravitation. The moon's orbit, under a close study, is  in fact complicated and very irregular, and the problem of  describing it is a typical instance of the so-called 3-body problem  on which several  mathematicians other than Clairaut were working. Among them we mention  Euler, d'Alembert, Maupertuis, Daniel Bernoulli and Gabriel Cramer.  In 1749, Clairaut was awarded a  prize a the competition organized by the Saint Petersburg Academy of Sciences for a memoir whose subject was the motion of the moon.

On February 6, 1744, Euler, who was the director of the class of mathematics at the Prussian Academy of Sciences, presented Clairaut's book \emph{Théorie de la figure de la terre} to this Academy, and at the same time he
 proposed the election of the author as a corresponding member.

Clairaut published in 1746 his treatise \emph{Elements of algebra} \cite{Clairaut-algebra}, which, like his \emph{Elements of geometry}, is written with his personal and unique style, following a logical pattern  but with no theorem or  problem stated.

 In 1751,  Clairaut won a prize set by the Saint Petersburg Academy of Sciences for a memoir titled    
   \emph{Théorie de la lune déduite du seul principe d'attraction réciproquement proportionnelle aux carrés des distances } (Theory of the moon deduced only  from the attraction principle inversely proportional to the squares of the distances). The subject of the competition was suggested by Euler, and the question set for it was: \emph{Are the discrepancies observed in the motion of the moon compatible with the Newtonian theory?}  Clairaut's winning-prize memoir was printed in 1752 in Saint Petersburg \cite{C1752}.

      At the same epoch, Euler obtained two times (in 1748 and 1752) the prize proposed by the Paris Royal Academy of Sciences. Clairaut was part of the jury for these prizes.  

   In 1750, Clairaut won a prize set by the \emph{Académie Royale des Sciences, Inscriptions et Belles-Lettres} of Toulouse, for a memoir titled \emph{Nouvelle théorie de la figure de la terre où l'on concilie les mesures actuelles avec les principles de la gravitation universelle} (New theory of the figure of the Earth where we conciliate the present measures with the principles of universal gravitation) \cite{Clairaut-Nouvelle}.

Among the other works related to our topic that Clairaut published in the next few years, we mention his famous \emph{Tables de la lune calculées suivant la théorie de la gravitation universelle} (Tables of the moon, calculated according to the theory of universal gravitation), first published in 1754 \cite{Tables1754}, and his  \emph{Mémoire sur l'orbite apparente du soleil autour de la Terre, en ayant
égard aux perturbations produites par les actions de la lune et des
planètes principales} (Memoir on the apparent orbit of the sun around the Earth, taking into consideration the perturbations produced by the actions of the moon and the main planets),  published in 1759  \cite{Clairaut-1759}.

   In 1754, Clairaut was elected corresponding member of the Academy of Sciences of Saint Petersburg.

Clairaut died in 1765, at age 52.     In the last ten years of his life, the papers he published were mostly on astronomy and optics, in particular, on the improvement  of telescopes.  During the same period, Euler was working on the same topics. The chair of Mechanics that Clairaut left  at the Academy of Sciences was occupied after his death by d'Alembert.

  \section{On Clairaut's works} \label{s:works}

In the previous section, we already mentioned several works of Clairaut. In particular, we mentioned his two memoirs on space curves,  \emph{Researches on the curves with double curvature}  \cite{C0}  which he wrote at age 16, and  \emph{On the curves that are formed by cutting an arbitrary curved surface by a position plane} \cite{C01} which he wrote a few years later. He wrote other articles on this subject, and we mention   among them
  \emph{Quatre problèmes sur de nouvelles courbes} (Four problems on new curves) \cite{C02}, published in 1734,     \emph{Des épicyclo\"\i des sphériques} (On spherical epicycloids)     \cite{C03}, published in 1735,  and  \emph{De la spirale d'Archimède décrite par un mouvement pareil à celui qui donne la cyclo\"\i de et sur quelques autres courbes de même genre} (On Archimedes' spiral described by a motion similar to the one which gives the cycloid and on some other curves of the same kind)
  \cite{C1742}, published in 1742.

We also mentioned Clairaut's theory of the motion of the moon. This theory, like the theory of the revolution of the planets around the sun and the rest of the science of astronomy in the eighteenth century, is based on Newton's law of universal attraction, the same one  that governs the falling bodies, and the one that was used by Newton, and after him, Clairaut and others, to explain the figure of the Earth. On this subject, and more generally on Clairaut's work on Newton's theory, we mention his memoirs \emph{De l'orbite de la lune dans ls système de M. Newton,} (On the orbit of the moon in the system of Mr. Newton) (1747) \cite{Clairaut-Lune} and  \emph{Du système du monde dans les principes de la gravitation universelle} (On the system of the world in  the principles of universal gravitation) (1749) \cite{Clairaut-systeme}.

  We shall mention now other works, with an emphasis on those concerning the figure of the Earth.

On December 5, 1733, Clairaut read to the Royal  Academy a memoir with a geographical title, \emph{Détermination géométrique de la perpendiculaire à la méridienne tracée par M. Cassini avec plusieurs méthodes d'en tirer la grandeur et la figure de la terre} (Geometrical determination of the perpendicular to the meridian drawn by Mr. Cassini with several methods of extracting its length and the figure of the Earth)  \cite{C1}. A few years later, he read another memoir, 
 \emph{Sur la nouvelle méthode de M. Cassini pour connaître la figure de la terre} (On the new method of Mr. Cassini to know the figure of the Earth)  \cite{C2}.
 Despite the titles, these memoirs are geometrical. In fact, motivated by the theory of the form of the Earth, Clairaut worked on  the geometry and trigonometry of the ellipsoid of revolution, instead of that of the sphere. His two memoirs \cite{C1} and \cite{C2}
contain several theorems on the geometry of the geodesics on a surface of revolution which is not the sphere, in particular on the curvature of the curves obtained by the intersections of such surfaces with planes, etc. with a special attention to the case of an ellipsoid of revolution. In the memoir \cite{C1} he obtains a trigonometric formula for a right triangle on such a surface.
 In 1734, he wrote a memoir titled \emph{Solution de plusieurs problèmes où il s'agit de trouver des courbes dont la propriété consiste
dans une certaine relation entre leurs branches, exprimée par une équation donnée} (Solution of several problems where it is asked to find curves whose property consists of a relation between their branches, expressed by a given equation)  in which he finds plane curves defined by certain geometrical conditions (see \cite{Clairaut1736}, published in 1736).
  Problem III is the following (see Figure \ref{MON}): 
  \begin{quote}
  \emph{To find the curve MON such that when we slide along it a right angle $MCN$, the vertex $C$ of this right angle stays on a given curve.}
\end{quote}

Christian Goldbach, in a letter to
 Euler dated October 12th, 1744 
informs the latter that he lately came across  this memoir of Clairaut and he says, concerning Problem III:  `` In my opinion no more general solution than the one he indicates can be
thought of, and his statements about Problem III also appear to be
very remarkable." \cite[Letter No. 84]{Euler-Goldbach}
%
%    
%\begin{figure}[htbp]
%
%\centering
%\includegraphics[width=10cm]{MON.pdf}
%\caption{\small{Figure for Problem III of Clairaut's \emph{Solution de plusieurs problèmes} \cite{Clairaut1736}.}} \label{MON}
%\end{figure}
%

Euler responded to Goldbach on  November 17th, 1744, saying: 
``At your instigation, Sir, I read Mr. Clairaut's paper in the 1734 \emph{Mémoires}.
The solution of the first two problems is the Newtonian one and could indeed
not be more general. The third problem is indeed very remarkable." \cite[Letter No. 85]{Euler-Goldbach}

     In his book, \emph{Théorie de la figure de la terre, tirée des principes de l'hydrostatique}  (Theory of the figure of the Earth, drawn from the principles of hydrostatics) \cite{Clairaut-Terre}, published in 1743, Clairaut studies the Earth as a spheriod,   based on the laws of hydrostatics, that is, the theory of the  equilibrium of forces acting on a fluid. He develops the idea that the Earth, originally constituted by a fluid matter, acquired gradually the form that is required by the equilibrium laws of hydrostatics.  In a certain sense, his theory of the figure of the Earth is an extension of Newton's theory of universal attraction. This book contains the so-called  Clairaut theorem in which the spheroid, representing the surface of the Earth, appears as the surface  in hydrostatic equilibrium, submitted to the sum of gravitational and centrifugal potentials satisfying a certain exact differential equation for a homogeneous field.  In this work, Clairaut, while he confirmed Newton's conclusions on the form of the Earth, corrected some of the latter's computations. 
     
 The first part of the work \cite{Clairaut-Terre}  is titled \emph{Principes généraux pour trouver les hypothèses dans lesquelles les fluides peuvent être en équilibre, et pour déterminer la figure de la Terre et des autres planètes, lorsque la loi de la pesanteur est donnée.} (General principles for finding the hypotheses under which fluids can be in equilibrium, and for determining the figure of the Earth and of the other planets, given the law of gravity).
 
In \S 1, Clairaut says that a  fluid mass cannot be in equilibrium unless the efforts of all parts contained in a channel of an arbitrary figure which we imagine as traversing the entire mass cancel each other.  This is the first time that a principle of fluid equilibrium is stated.
  
  Clairaut corrected Newton's computations, showing that with the Earth considered as a fluid in equilibrium, the ratio of the minor axis to the major axis of the spheroid that it forms is 230/231, which is different from the value that Newton found.
From the mathematical point of view, the vector representing  gravity at an arbitrary point, that is, the attraction force which makes bodies fall, in the case of a spheroid, is not directed towards the center, unless the spheroid is a sphere.

In a note, on p. xiii of the introduction of this book, Clairaut declares that, following Maupertuis, he makes a distinction between the notions of  \emph{weight} (pesanteur) and gravity (gravité).
The word weight designates  the natural force that makes an arbitrary body fall, whereas
 gravity is the force with which this body falls if its effort and direction were not altered by the Earth rotation.
A famous result of Clairaut gives the value of this gravity (p. 249-250 of the book).  

The methods that Clairaut uses in the work are geometrical. A few years later, with Lagrange's \emph{Mécanique analytique} (1788), analytical methods became dominant in mechanics.

 The Earth as a spheriod was sometimes referred to in the eighteenth century science literature as the \emph{Clairaut spheroid}.
 For an exposition of Clairaut's ideas in relation with those of Newton, we refer the reader to the survey  \cite{Taton-Diffusion} by René Taton.

      The memoir
 \emph{An Inquiry concerning the figure of such planets as revolve about an axis, supposing the density continually to vary, from the centre towards the surface}   \cite{Inquiry} was presented to the Royal Academy of Paris and to the Royal Society of London in 1737-1738, and it was published in the latter's Transactions, translated into English. In this memoir, Clairaut returns to Newton's part of the \emph{Principia} in which the latter is concerned with the figure of the Earth. He first announces that some observations he made under the Arctic circle led him to believe that this figure was flatter than what Newton thought. He then expresses his surprise concerning the fact that Newton applied different physical theories, as to the causes of this ellipticity, regarding the Earth and  Jupiter. But Clairaut is mostly interested in  geometry, and the core of his memoir is mathematical. In fact, his aim is to discusses the geometrical problems that are motivated by those of the figure of the Earth. He declares (p. 179):
 \begin{quote}\small
And though my hypothesis should not be conformable to the laws of nature, or even though it should be of no real use [\ldots], I thought however that geometers would be pleased with the speculations contained in this paper, as being, if not useful, yet curious problems at least.
 \end{quote}

 Among the problems that Clairaut discusses, we quote the following three:

   \begin{quote}\small
    \emph{Problem 1: To find the attraction which a homogeneous spheroid, differing but very little from a sphere, exerts upon a corpuscle placed at a point on the axis of revolution.}
    
    \medskip

    \emph{Problem 2: The spheroid is no more supposed to be of a homogeneous matter, but composed of an infinite number of ellipsoidal strata which are all similar, and whose densities are represented by an arbitrary curve of which the equation is known. To find the attraction that it exerts on a corpuscle placed at a pole.}
    \medskip

          \emph{Problem 3: To find the attraction which a spheroid exerts upon a corpuscle placed at an arbitrary point of its surface.}
          
          \end{quote}
          
       Clairaut's work is based on Netwon's theory of universal gravitation, but also on  Huygens' work on the centrifugal force.  He writes: ``[Huygens' principle] which consists in making bodies gravitate perpendicularly to the surface, seems to me of absolute necessity".

                    Among the other memoirs of Clairaut, we mention his 
 \emph{Détermination géométrique de la perpendiculaire à la méridienne} (Geometrical determination of the perpendicular to the meridian) \cite{C1} in which he shows that such a perpendicular is never a plane curve, unless the Earth is considered as spherical. In the other cases, any perpendicular to a meridian, except the equator, is a curve with double curvature, that is, it is not contained in a plane. He also studies properties of this curve. 
      
     From the point of view of fluid mechanics, we quote Lagrange, who declares in his \emph{Mécanique analytique} that Clairaut changed the face of Hydrostatics, and made it a new science. He writes \cite[t. 1, p. 179-180]{Lagrange-Mecanique}:

    \begin{quote}\small
     [\ldots]  Clairaut made [Newton's principle] more general, by showing that the equilibrium of a fluid mass requires that the forces of all the parts of the fluid enclosed in an arbitrary channel, ending at the surface or entering into itself, destroy each other. He was the first to deduce, from this principle, the true fundamental laws of equilibrium of a fluid mass whose parts are animated by arbitrary forces, and he found the partial difference equations by which these laws can be expressed, a discovery which changed the face of Hydrostatics, which he made a new science.
\end{quote}   
    
    After this passage, Lagrange talks about Euler's work on hydrostatics, which, he says, is adopted in almost all the treatises on this science.
    We shall discuss some of Euler's works related to the figure of the Earth, but before that, we quote d'Alembert on this question.

  \section{The figure of the Earth in the d'Alembert--Diderot encyclopaedia} \label{s:enc}
  The sixth volume of Diderot and d'Alembert's  \emph{Encyclopédie}, published in 1756, contains a  long article on the figure of the Earth. It is written by d'Alembert.  I would like to quote a few passages (which I am translating into English), which show the importance of this question in the eighteenth century, and in which Clairaut is described as one of the main figures in this debate. D'Alembert also mentions the relation between the latter's work and that of the mathematician Colin Maclaurin, a question which is also addressed in the correspondence between Euler and Clairaut. 
\begin{quote}\small 

 This important question has made so much noise in recent times and the scholars, especially in France, have been so engaged in it, that we thought we should devote to it a special article.
  
[\ldots] 
Mr. Clairaut having meditated on this last condition, deduced from it
consequences, which he presented in 1742 in his treatise titled \emph{Theory of the figure of the 
Earth, derived from the principles of Hydrostatics}. According to Mr. Clairaut, in order for a fluid to be in equilibrium, the forces of all the parts included in a channel of any shape that one can imagine crossing the whole mass, must cancel out each other. 
This principle is apparently more general than that of Mr. Maclaurin; but I have shown in my essay on fluid resistance  (1752. art. 18) that the equilibrium of curvilinear channels is only a corollary of the simpler principle of the equilibrium of rectilinear channels of Mr. Maclaurin. This however does not diminish anything of the merit of Mr. Clairaut, since the latter deduced from this principle a great number of important truths that Mr. Maclaurin had not drawn from it, and of which he had even been enough unaware  to fall into some errors; for example, in those of supposing the layers of a spheroidal  fluid  to be similar to each other, as can be seen in his treatise on fluxes, art. 670. and seq.

Mr. Clairaut, in the work we have just quoted, proves (a fact which Mr. Maclaurin did not do directly) that there are  infinitely many hypotheses in which the fluid would not be in equilibrium, even though the central columns would counterbalance each other, and even though gravity would be perpendicular to the surface. He gives a method for recognizing the hypotheses of gravity in which a fluid mass can be in equilibrium, and for determining its figure; he proves moreover that in the system of attraction of the parts, provided that gravity is perpendicular to the surface, all the points of the spheroid will be equally pressed in any direction, and therefore that the equilibrium of the spheroid, in the hypothesis of attraction, is reduced to the simple law of perpendicularity to the surface. According to this principle, he seeks the laws of the figure of the Earth under the hypothesis that the parts attract each other, and that it is composed of heterogeneous layers, either solid or fluid;
he finds that the Earth must have in all these cases a more or less ``flattened" elliptical figure, according to the disposition and the density of the layers: he proves that the layers must not be similar, if they are fluid; that the increases of gravity from the equator to the pole must be proportional to the squares of the sines of latitude, as in the homogeneous spheroid; this is a very remarkable and very useful proposition in the theory of the Earth: he proves moreover that the Earth could not be more flattened than in the case of homogeneity, that is to say, from 1230; but this proposition only takes place by assuming that the layers of the Earth, if it is not homogeneous, increase in density from the circumference towards the center; a condition which is not absolutely necessary, especially if the inner layers are assumed to be solid; moreover, even supposing that the densest layers are the ones closest to the center, the flattening may be greater than 1230, if the Earth has a solid inner core more flattened than 1230. See the third part of my \emph{Recherches sur le système du monde}, p. 187. Finally, Mr. Clairaut demonstrates, by a very beautiful theorem, that the decrease of gravity from the equator to the pole is equal to twice 1230 (the flatness of the homogeneous Earth) minus the real flatness of the Earth. This is only a very slight sketch of what can be found as excellent and remarkable in this work, which is very superior to all that had been done until then on the same subject.

\end{quote}
%   \begin{figure}
%\centering
%\includegraphics[width=11cm]{Layers.pdf}
%\caption{\small{The stratification of the Earth into different layers of various densities. Figure extracted from Clairaut's  \emph{Theory of the figure of the Earth} \cite{Clairaut-Terre}.}} \label{fig:stratification}
%\end{figure} 

  \section{On Euler's works related to the figure of the Earth}\label{s:Euler}

 The issue of the figure of the Earth, at the age of Enlightenment, was discussed in the literary and learned circles and not only in the academic milieux. Indeed, the problem was easy to understand by the laymen. Furthermore, it was addressed in several pamphlets of Voltaire, the preeminent philosophical and literary figure of the Enlightenment who was also one of the heaviest French proponents of Newton's theories, as opposed to those of Descartes. Voltaire also co-authored a book titled  \emph{Elements of the Philosophy of Newton}  with the philosopher \'Emilie du Châtelet who was at the same time his collaborator and mistress.  The book played a decisive role  in the popularization of Newton's ideas in France.
\'Emilie du Châtelet  also published a French translation with commentaries of Newton's \emph{Philosophiae Naturalis Principia Mathematica}.  Voltaire wrote the preface to this work. We read in this preface, concerning du Châtelet's commentary:
\begin{quote}\small
Regarding the Algebraic Commentary, it is a work beyond the translation. Madame du Châtelet worked on it following the ideas of Mr. Clairaut: she made all the calculations herself, and whenever she had finished a chapter, Mr. Clairaut examined it and corrected it. Not only that, it is possible for a mistake to be made in such a difficult work: it is very easy to substitute one sign for another when writing; Mr. Clairaut also had a third person review the calculations when they were finalized, so that it is morally impossible for an inattention error to have slipped into this work; and it would be just as impossible, for a work in which Mr. Clairaut participated, not to be excellent in its kind.
 \end{quote}
   
 At about the same period, Euler published a paper also addressed to the general public, in 7 installments, titled   \emph{Von der Gestalt der Erden}   (On the shape of the Earth),  in a magazine published in German in Saint Petersburg, the  \emph{Anmerckungen \"uber die Zeitungen}  (Notes on the newspapers).  In this paper,  Euler surveys the various problems that were addressed at that time by the question of the figure of the Earth, with  the three conflicting points of view: the Earth is either spherical, or spheroidal and flattened at the poles (and he says that in this case it has the shape of an orange), or elongated at the poles (and he says in this case that it has the shape of a melon).

 Euler recalls in this paper that in the seventeenth century, supporters of each of the two latter views had strong arguments, and he surveys the ideas from physics that lie behind these hypotheses. 
  He mentions that, mathematically, one has to work  under the assumption that the surface of the Earth is smooth, for instance, that it consists of still water. Indeed, to use geometry and differential calculus, one needs to introduce at each point vectors that are perpendicular to the surface of the Earth, representing the action of gravity. If  the Earth is melon-shaped, a body would be heaviest if it is closer to the Equator than to the poles, and  lighter in the inverse situation. He notes that the various opinions that were made concerning this question were not only based on experiments and measures, but also on deep theories and he recalls that this question is important not only for geography but also for physics. 
He mentions the experiments with pendulums that were done in relation with this question. In the case of a melon-shaped figure of the Earth, and assuming these pendulums beat accurately a second, they should have a greater length at the equator than at the poles. He describes in detail the geometry that is behind the measurements of the degrees of the meridians that were conducted near the poles and at the equator with the help of astronomical observations. He says that even if, at the time he was writing his article, the results of the latter expedition were still not available, he is quite confident they they will confirm the fact that the Earth is orange-shaped. He says that the measurements done during of the Peru expedition will still be important  because they will give information on the length of the diameter of the Earth.

  Let us make a quick review of some other writings of Euler on this subject.

The title of the following memoir is quite informative:  \emph{Methodus viri celeberrimi Leonhardi Euleri determinandi gradus meridiani pariter ac paralleli telluris, secundum mensuram a celeb. de Maupertuis cum sociis institutam} (Method of the celebrated Leonhard Euler for the determination of a degree of a meridian, as well as of a parallel of the Earth, based on the measurement undertaken by the celebrated de Maupertuis and his colleagues) \cite{E132} (1750). In this memoir, Euler considers several geographical problems, under the assumption that the Earth is a spheroid, including the question of  determining the degree of a meridian given the elevation of the pole,\footnote{Euler sometimes used the term ``elevation of the pole" to denote altitude. The term originates from astronomy. See his memoir \cite{E215}, translated in \cite{CP}.} and the one of  finding the elevation of the pole for a given degree of a parallel.  At the same time, he reviews some work of Christian Nikolaus Winsheim,\footnote{Christian Nikolaus von Winsheim (1694-1751) was an astronomer and geographer who settled in Saint Petersburg in 1718. In 1731, he became adjunct in astronomy  at the Imperial Academy of Sciences and  associate dean of the astronomy section. He was in charge of the calculations that were necessary for the observations. It was not possible for him to participate in these observations since, because of his obesity, he was not able to climb  the steep spiral staircase which led to the 4th, 5th and 6th floors of the tower of the \emph{Kunstkamera} in Saint Petersburg, where the observatory was located. 
He collaborated with Euler on the preparation of the \emph{Russian Atlas} published in 1745 \cite{Atlas-Russicus}. In 1745, he published a book considered as the first textbook on political geography.} who was his colleague at the geography section of the Saint Petersburg Academy of Sciences.  In the appendix to this memoir, Euler  provides tables of lengths of degrees of  meridians and parallels at various places.
 
   The next memoir we mention is connected with the influence of the figure of the Earth on astronomy. It was published in 1747,  and it is titled \emph{De la parallaxe de la lune tant par rapport à sa hauteur qu'à son azimuth, dans l'hypothèse de la terre sphéroïdique} (On the parallax of the moon with respect to its elevation and azimuth, under the hypothesis of a spheroidal Earth) \cite{E172}. The term ``parallax" refers to the  influence of the position of an observer on the trajectory of a celestial object, seen from his own position (the object being, here, the moon).  Note that the notion of parallax applies only to the observation of objects within the solar system, otherwise, the observation does not depend on the position of the observer on Earth, since the Earth can be assimilated to a point. From the mathematical viewpoint, the problem studied in this memoir is a coordinate change problem, in the new spheroidal geometry setting. Euler starts his memoir by recalling that Maupertuis published an ``excellent treatise" on the parallax of the moon, in which he showed how the usual rules (under the hypothesis of a spherical Earth) have to be modified, if one takes into consideration the spheroidal shape, but that the latter failed in taking onto account one  parameter, namely the azimuth, that is, the angle seen from the observer, in a horizontal plane, between the projection of the direction of the celestial object considered, and a given reference direction. He then develops the trigonometrical computations needed in this geometrical problem.  
   
In the memoir \emph{De attractione corporum sphaeroidico-ellipticorum}
 (On the attraction of spherico-elliptical bodies) \cite{E97}, published in 1738,  Euler gives a formula for the law of attraction between a particle  situated at a pole and another one situated at the equator, in the case of a planet made of a uniform material, where the particles attract each other according to Newton's law, that is,  by a force whose magnitude is inversely proportional to the squares of the distances and which rotates about the axis.

We then mention Euler's important memoir   \emph{\'Elémens de la trigonométrie sphéroïdique tirés de la méthode des plus grands et plus petits} (Elements of spheroidal trigonometry drawn from the method of the maxima and minima) \cite{E215}.  Euler discusses there the mathematical theory that is behind the measures of degrees of meridians that were conducted during the Peru and Lapland expeditions, the possible errors made during these measurements and their impact on the knowledge of the true figure of the Earth. At the same time, he develops the geometry of a spheroid. Figure \ref{fig:Elements1} is extracted from that memoir.
Among the mathematical results contained in this memoir we highlight the following:

   \begin{enumerate}

\item Given the latitude of a point, to determine its distance to the centre of the Earth  (Section 5 of  \cite{E215}).

    \item  \emph{Given two points situated on the same meridian and whose latitudes are known, to find their distance} (Section 17 of  \cite{E215}).

\item \emph{Given two points $L$ and $M$ of which we know the latitudes and the difference of their longitudes, to find the shortest path $LM$ on the surface of the Earth which leads from one point to the other}  (Section 19 of  \cite{E215}).

\item   \emph{To determine the ratio of the diameter of the Equator to the axis of the Earth, without using the measures done by the expeditions near the pole and near the equator, but by a construction done in a small portion of land}
 (Section 24 ff. of  \cite{E215}).

 \end{enumerate}
 
 For the last problem, Euler proposes a series of astronomical observations, together with the possibility of drawing a straight line (a geodesic) in the given region. Provided this can be done precisely, Euler gives a formula for the value of the ratio of the diameter of the Earth to its axis. 
 
 Note that the answers to the first three problems are straightforward in the case of a sphere. In the case of a spheroid, the latitude  has to be defined with care; this is the angle made by a perpendicular to the surface of the Earth with its axis of rotation, at the intersection point of these lines;  note also that this perpendicular does not pass by the center of the spheriod, unless the point is on the equator. 
 See Figure \ref{fig:Elements1} in which this angle is denoted by $ENM$. Note also that this perpendicular is the direction of the gravity.
 
% \begin{figure}[htbp]
%\begin{center}\label{fig:Elements1}
%\input{Euler_Elements-fig1.pdf_tex}
%\caption{$AB$ is the axis of the Earth and $CE$ its semi-diameter.}
%\end{center}
%\end{figure}
%

  In Section 12 of the memoir, Euler, after a computation, gives an estimate of the ratio of the diameter of the Earth to its axis, which is  230/229, and he notes that this value coincides with the one found by Newton. He writes:  ``Whence we can conclude that the hypotheses  made by this great geometer on the structure and the attraction of the Earth agree with reality." In the sections that follow, he discusses other estimates obtained by using other methods of computation, and results of measurements of degrees of meridians conducted  at various places on the surface of the Earth. The methods he gives are successful in regions that are not too close to the equator or to the poles (\S 34).

    To end this section, I would like to mention three memoirs of Euler on hydrodynamics, a topic in which Euler was greatly motivated by his discussions with Clairaut on the formation and the figure of the Earth.
    The three memoirs are considered as having completely transformed
    the field of fluid dynamics.
         \begin{enumerate}
    \item Principes généraux de l'état d'équilibre des fluides  (General principles of the state of equilibrium of fluids) \cite{E225};

      \item   Principes généraux du mouvement des fluides (General principles of the motion of fluids) \cite{E226};
      
\item
Continuation des recherches sur la théorie du mouvement des fluides
(Continuation of the researches on the theory of the motion of fluids) \cite{E227}.
    \end{enumerate}

  In the first memoir, Euler makes a thorough study of the notion of pressure and its applications in fluids, based on a heavy use of differential calculus and mathematical analysis.  He  presents for the first time the general equations of hydrostatics and the equilibrium equations for fluids. In the  second memoir, Euler develops several theories emitted in the first, analyzing  for the first time vortex flows and Poiseuille flows.\footnote{Poiseuille flows are  also known as Hagen--Poiseuille flows. They are  named so after the 19th century  French
  physicist Jean-Léonard-Marie Poiseuille and the Prussian engineer Gotthilf Hagen. These are flows that follow the so-called Poiseuille law which makes a relation between the amount of flow of a fluid, the viscosity of the fluid in a pipe, the difference in pressure at the boundary of the  pipe, and the length and the diameter of this  pipe.} The third memoir is a continuation of the first two, in which Euler improves one his previous results, treating  compressible fluids and developing an exhaustive theory of flows in pipes.  The three memoirs are thoroughly analyzed by C. Truesdell in his  introduction to the volume of Euler's \emph{Opera Omnia} (Ser. II, Vol. 12) in which these three papers are reproduced.

   \section{On the Euler--Clairaut correspondence} \label{s:correspondence-C}

The Euler--Clairaut correspondence published in Euler's  \emph{Opera Omnia} contains 60 letters,  among which 15 were sent by Euler to Clairaut and 45 by Clairaut to Euler. The topics discussed include integration, isoperimetry, differential calculus, differential geometry of curves and surfaces, functions of several variables, differentiation of integrals depending on a parameter, number theory,   Newton's theory of attraction, hydrostatics, mechanics, elasticity,  engineering problems,  celestial mechanics (in particular  the motion, of the moon and of  Saturn),  questions related to the figure of the Earth, and other topics.

When, on  September  17, 1740, Clairaut sent his first letter to Euler, he was 27 years old and Euler was 34. This was the year before Euler left Saint Petersburg to Berlin. (Euler had arrived to Saint Petersburg at age 20). 
Clairaut was already a well-known mathematician in Paris, who had already published on several subjects that were among Euler's interests. At that time, Clairaut had already close relations with Daniel Bernoulli, who was a close friend of Euler, and with Pierre-Louis de Maupertuis, who was also in contact with Euler. 
Euler was aware of Clairaut's work;  Daniel Bernoulli, in a set of letters sent to Euler  in 1739, had mentioned several works of Clairaut  in a very laudatory way, and Euler  had asked Bernoulli several details about these works, see \cite[note 3 p. 69]{Opera-4-A-7}.

 The first letter from Clairaut to Euler is dated 17 sept 1740. We already quoted the beginning of this letter. The question of the figure of the Earth is not mentioned in this letter. Clairaut discusses integration and differentiation, and in particular, he asks Euler questions concerning integrals of the form
\[\int \frac{xadx}{x^2+a^2+ax}
\]
in which the variable is $a$.

 I will now review some excerpts of the letters that are related to the topic  of the figure of the Earth. The discussion on this question starts  in Euler's first letter to Clairaut.

Euler responded to Clairaut's first letter on 19 octobre 1740. Euler writes that he himself had made a great progress on this matter thanks to the reflections of Clairaut which had reached him, and he starts by mentioning the memoir that the latter wrote on the figure of the Earth, published in the \emph{Transactions of the Royal Society} \cite{Inquiry}. He tells his correspondent that these results are not only interesting in themselves, but that they will give us information on the internal structure of the Earth. He informs him that Maupertuis found recently that the ratio of the axis of axis of the Earth to the diameter of the equator is 177/178, i.e., it is greater than the one which would follow from the hypothesis that the Earth is homogeneous, and  that this implies that the density of the Earth is greater on the surface than at the center, which is contrary to Newton's opinion expressed in Book II, Prop. X of the \emph{Principia}. 
 
 In his response, dated December 26, 1740, Clairaut tells Euler, concerning this matter, that he trembled at the idea that he might have made a mistake in issuing a statement contrary to Newton's opinion, especially because he learned from (Daniel) Bernoulli that the same statement was also contrary to Euler's opinion. We note that later, in his \emph{Théorie de la figure de la terre} (1743), Clairaut found that the flattening of the Earth at the poles should be around 1/230.\footnote{It turned out later (after the precise measures done in the XXth c.) that the value is approximately 1/298, see \cite[p. 77, Note 3]{Opera-4-A-7}.} 
 
 In the same letter, Clairaut talks about the competition set by the French academy, to which Euler had participated, and whose subject was the form the ebb and flow of the sea.  Clairaut was part of the jury. The prize went to Euler, A. Cavalleri, D. Bernoulli and C. Maclaurin. The four artciles were published by the French Academy. Clairaut tells  Euler that Mauclaurin has a beautiful proof of the fact that Earth has the form of an ``Apollonian ellipse", under the assumption that all its parts are under a mutual attraction which is inversely proportional to the square of the distance.

On  8 January 4, 1742,  Clairaut announced to Euler  that he finished writing his book \emph{Théorie de la figure de la terre} \cite{Clairaut-Terre}, insisting on the physical part, that is, on fluid mechanics. He tells him about corrections that have to be done on the paper he published on the figure of the Earth in the Philosophical Transactions \cite{Inquiry}, in which he assumed that the Earth was made of successive elliptical layers, each having its proper density. We note incidentally that Huygens, when he introduced in the problem of the figure of the Earth the centrifugal force due to the motion of the Earth around its axis, considered that the Earth density is uniform. Clairaut thought that these layers are not similar, but that their flatness increases with their distance to the center of the Earth. The layers are submitted to the laws of attraction among each other and to the motion of the Earth around its axis. Starting from this letter, Clairaut  informed regularly Euler of his work on the subject, from the point of view of mechanics, in particular his theory of the attraction field in a spheroid.

In a letter sent from Euler to Clairaut in January-February 1742 (no precise date), Euler congratulates the latter for his book \emph{Théorie de la figure de la Terre}, which, he tells him, ``will enlighten as well the so important question of the figure of the Earth, that it will increase the admiration that all the world conceived for you". Euler writes in particular: 
 
\begin{quote}\small
Your remarks on the variation of the types of ellipses that form the different layers of the Earth are extremely profound, especially if they are not a consequence of a hypothesis, but of the indisputable principles of hydrostatics, as you have the goodness to assure me. Moreover, your piece in the Transactions  still has its merits, because up to now it seems to me that it has been possible to form hypotheses based on it which are suitable for calculations.
\end{quote}

In his response, sent on March  28, 1742, Clairaut gives Euler more details on the problem of equilibrium of forces for fluids. At the same time, he writes: ``I am very flattered by the desire you show to read my work on the figure of the Earth, but I fear that it does not respond to the advantageous idea you have of it."

He explains that he has two  principles for determining the gravity assumptions that are possible in fluids. The first one says that a channel of an arbitrary shape going from one point of the surface to another is in equilibrium independently of the rest of the fluid. The second one says that if we assume that the fluid mass is completely divided into infinitely many layers by surfaces at each point of which gravity is perpendicular to the tangent plane, then at an arbitrary point $M$ on any one of these layers, the thickness $MN$ is inversely proportional to the force of gravity at the point $N$. 
From these two principles he deduces the following 
 
\begin{quote}\small
\noindent \emph{Theorem.--- Suppose that the force of gravity is decomposed into two others, the one parallel to the axis and the other to the equator. Furthermore, let the first force be expressed by a function $P$ of the coordinates $CQ$, $x$ and $QM$, $y$ and the second by another function $Q$ of the same variables. Then,  $Pdy+Qdx$ is the differential of some function of $x$ and $y$ for the fluid to be in equilibrium. In other words, we necessarily have $dP/dx=dQ/dy$.}
\end{quote}
Clairaut tells Euler that if he were not tired, he would have explained to him the proof of this theorem. The proof appeared in Clairaut's \emph{Théorie de la Figure de la Terre}.  
 
%
%\begin{figure}[htbp]
%
%\centering
%\includegraphics[width=11cm]{letter.pdf}
%\caption{\small{Figure from Clairaut's letter to Euler, dated March  28, 1742}} \label{fig:letter}
%\end{figure} 
The discussion on hydrodynamics, motivated by the question of the figure of the Earth, continues in the  letter from Euler to Clairaut, sent in April, 1742, and the one from Clairaut to Euler, sent on May 29 1742, in which incidentally Clairaut thanks Euler for his worries about his health, and he tells him that contrary to what Euler suspects, it is unlikely that these health problems are due to an excess of work. In this letter, Clairaut expands on questions of fluid dynamics and on the problem of the figure of the Earth.

 In his letter  dated July 25, 1742, Clairaut writes to Euler that nothing that has to do with the nature of fluids is clear, because in this field one cannot have demonstrations as in Geometry. He tells him that, because of this fact, he fears that the book he intends to  publish on this subject will not fully satisfy him but that in any case, he will receive with pleasure his comments on it.
 
  In the  letter dated December 3, 1742, Clairaut  tells Euler that he is spending most of his time finishing his book on the figure of the Earth. He also writes that he has read a chapter written  by Maclaurin on the same subject. He compares at length his results with those of Maclaurin and he declares that after having examined his writings, he considers that the latter's  assumptions on fluid dynamics are not correct.
 
 In his letter dated April 23, 1743,  Clairaut announces to Euler that his book on the Figure of the Earth is published, and that he is looking forward for his reaction on it.

Euler and Clairaut  never met. This was not unusual in those times. Likewise, Euler and Lagrange, although they had a large variety of common interests, never met.   

 Among the other mathematicians with whom  Clairaut had a regular correspondence, we mention  Johann I and Daniel Bernoulli,  Gabriel Cramer and Colin Maclaurin.  The correspondence between scientists is an invaluable source of information on their work and on the work of other mathematicians, especially in those times, where there were still very few mathematical journals.

The magazine  \emph{Leipzig Zeitungen von gelehrten Sachen} published a review of  
Clairaut's book \emph{Théorie de la Figure de la Terre} in its issue of
May 4th, 1744 \cite{Review1744}. Christian Goldbach, in a letter to Euler sent from Moscow, dated 
June 1st, 1744, \cite[p. 838]{Euler-Goldbach}, writes: ``From the review which the \emph{Leipzig Zeitungen von gelehrten Sachen} gives of
the \emph{Théorie de la figure de la Terre} by Mr. Clairaut, I gather it has to be a very
good book." Euler responded to Goldbach on  July 4th of the same year, saying: ``Mr. Clairaut's \emph{Théorie de la figure de la Terre} is indeed an incomparable
work, both with respect to the profound and difficult questions treated in it and
to the pleasant and easy method by which he is able to present the most sublime
matters very plainly and clearly."

  \section{Later works on the geometry of the spheroid}\label{s:works}

It is well known that some of the most preeminent mathematicians were also geographers. We mentioned Euler and Clairaut, but there are many others. We recall that 
 C. F. Gauss was the head of the observatory of the University of G\"ottingen and as such he had also the official title of geographer.  In 1820,
 George IV, who had the title of King of the United Kingdom of Great Britain and Ireland and also King of Hanover, gave him the task of 
 measuring  the extent of the German kingdom of
Hanover, to which the city of G\"ottingen belonged. Gauss, for his measurements, submitted this land to a triangulation, and used tastronomical observations for the determination of coordinates. 
In his paper \cite{Gauss- Allgemeine}, published in 1825, titled  
\emph{Allgemeine Auflösung der Aufgabe: die Theile einer gegebnen
Fläche auf einer andern gegebnen Fläche so abzubilden, daß die Abbildung dem
Abgebildeten in den kleinsten Theilen ähnlich wird} (General solution of the problem:
to represent the parts of a given surface on another so that the smallest parts of the
representation shall be similar to the corresponding parts of the surface represented), he shows that every sufficiently small neighborhood
of a point on an arbitrary real-analytic surface can be mapped conformally onto
a subset of the plane. In the same paper, which contains crucial mathematical results, he writes that his aim  is only to construct geographical maps and to study the
general principles of geodesy for the task of land surveying. Surveying the kingdom
of Hanover took nearly two decades to be completed. It led Gauss gradually to the
investigation of triangulations,  the use of the method of least squares in geodesy, 
and then to his major work, the \emph{Disquisitiones generales circa superficies curvas} (General investigations on curved surfaces).\footnote{It is in this article that Gauss proves the result he calls ``remarkable theorem" (\emph{Theorema Egregium}), which explains in particular why curvature is the only obstruction for a surface to be faithfully represented on the plane (\S 12; p. 20 of the English translation).} In the latter, we
can read,  in §27 (p. 43 of the English translation \cite{Gauss-Disquisitiones}): ``[\ldots] Thus, e.g., in the greatest of the triangles which we have measured in recent years, namely that between the points Hohehagen, Brocken, Inselberg, where the excess of the sum of the angles was 14".85348, the calculation gave the following reductions to be applied to angles: 
Hohehagen: -4".95113;
Brocken: -4".95104;
Inselberg: -4".95131."

F. W. Bessel, in 1837, wrote a paper titled \emph{Bestimmung der Axen des elliptischen Rotationssphäroids, welches den vorhandenen Messungen von Meridianbögen der Erde am meisten entspricht} (Determination of the axes of the elliptical rotational spheroid that is most consistent with existing measurements of Earth meridian arcs).

 In a paper published in 1841 and titled \emph{De la ligne géodésique sur un ellipsoïde, et des différents usages d'une transformation analytique remarquable} (On the geodesic line on an ellipsoid and the various usages of a remarkable analytic transformation) \cite{Jacobi}, C. G. J. Jacobi, studies geodesics on ellipsoids. These surfaces are more general than spheroids. He declares that his motivation for the study of this problem arises from geography, and he mentions works of Lambert and Gauss on this topic, as well as works by Euler on mechanics. Jacobi was led in this study to abelian integrals, which is one of his favorite subjects. 
It is also interesting to know that Jacobi studied similar problems of
geodesy using elliptic functions. In another paper published in 1857 and titled \emph{Solution nouvelle d'un problème
fondamental de géodésie} (A new solution of a fundamental problem in geodesy) \cite{Jacobi-Solution},
Jacobi considers, on an ellipsoid having the shape of the Earth, a geodesic arc whose
length, the latitude of its origin and its azimuth angle at that point are known, and he studies the
question of finding the latitude and the azimuth angle of the extremity of this arc,
as well as the difference in longitudes between the origin and the extremity. He writes: ``The problem of which I just gave a new solution has been recently the
subject of a particular care from Mr. Gauss, who treated it in various memoirs and
gave different solutions of it."

From the purely geometrical point of view, studying the differential geometry and the trigonometry of a  spheroid became gradually a fashionable subject, see \cite{Cayley, Forsyth1, Forsyth2, Gauss, Gundermann,  Helmer,  Ivory, Jacobi, Kruger, Kummell, Laplace-Oeuvres2, Legendre, Legendre-analyse, Monge, Poincare1905, Puissant1, Puissant2, Stein, Ward}.
We briefly mention a few results proved in some of these papers.

Laplace, in his \emph{Mécanique céleste} (Celestial mechanics) \cite[p. 128ff]{Laplace-Oeuvres2} studies on a spheroid the curves whose length is the shortest  distances between their endpoints.  
Legendre, in his \emph{Mémoire sur les opérations trigonométriques, dont les résultats dépendent de la figure de la terre} (Memoir on the trigonometric operations whose results depend on the figure of the Earth) (1787) \cite{Legendre} solves several questions on the geometry of the spheroid, one of them (\S VII-X of his memoir) concerning triangles and distances on the spheroid. In particular, he obtains simple formulae for the shortest line which starts at a given point, making with the meridian a given angle.  He writes that this work is motivated by problems in geography, in particular, those of establishing precise 
measurements for the coast of France.
He  gives several propositions concerning chains of triangles and triangulations. In a sequel to this memoir that he published several years later, titled \emph{Analyse des triangles tracés sur la surface d'un sphéroïde} (analysis of triangles drawn on the surface of a spheroid) \cite{Legendre-analyse} (1806),  he continues this study. He says that his work was motivated by the recent geodetical triangulation that were drawn to measure the distance between the cities of Dunkerque and Barcelona.

 In his paper \emph{Sur les lignes géodésiques des surfaces convexes} (On the geodesic lines of convex surfaces) \cite{Poincare1905}, Poincaré has a section (\S 2, p. 244-250) on the geodesics of the spheroid, which he studies in the setting of differential equations, using  a method that Lagrange  introduced in his study of the motion of a planet under the action of the perturbations due to the other planets, and which he calls the ``theory of variation of the constants." In some sense, this study is a generalization of a study that Poincaré made earlier of  geodesics on the sphere.\footnote{One must be careful here about the terminology: Here, Poincaré calls \emph{spheroid} a convex surface which is sufficiently close to a sphere. Thus, a spheroid in the sense of the present paper is a special case of a spheroid in sense of Poincaré.} The main questions in which he is interested is the number of stable closed geodesics on such a surface, the word ``stable" meaning here that the closed geodesic remains so under small deformations. In the case of the spheroid he is studying, he shows that this number is odd. After treating this question on a  spheroid,  he studies the same question in the case of a general convex surface. Poincaré was interested in these questions because of the relation (that he highlighted) with the 3-body problem.

In the more recent paper by Ward \cite{Ward}, the following question is studied: 

\begin{quote}
\emph{Given two points $P_1$ and $P_2$ on a spheroid which are not on the same meridian, let $s_{12}$ be the length of the geodesic connecting them, and $\sigma_{12}$ be the length of arc which is the intersection of the spheroid with the plane containing $P_1$ and $P_2$ and which passes by the center of the spheroid. By how much does $\sigma_{12}$ exceed $s_{12}$? }
\end{quote}

Ward declares at the beginning of his paper that the question is motivated by current interest in navigation.

 It is interesting to see how questions in geography contributed  to the development of geometry.

\end{document}